\documentclass[12pt]{article}
\usepackage{amsmath,amssymb,amsthm,eucal, mathtools, mathrsfs}
\usepackage[T1]{fontenc}
\usepackage{color}
\usepackage[all]{xy}
\usepackage[pdftex, draft=false]{hyperref} 
\title{\vspace*{-1pc}%
The intersection of complemented submodules is not complemented\\ \medskip
\small Addendum to: The Friedrichs angle and alternating projections in Hilbert $C^{*}$-modules}
\author{Koen van den Dungen{\dag\dag}, Bram Mesland\S$^*$, Adam Rennie\dag
\thanks{email: \texttt{kdungen@math.uni-bonn.de},
\texttt{b.mesland@math.leidenuniv.nl}, \texttt{renniea@uow.edu.au}
}
\\[3pt]
\dag\dag Mathematical Institute of the University of Bonn,  Germany
\\[3pt]
\S Mathematisch Instituut, Universiteit Leiden, Netherlands
\\[3pt]
\dag School of Mathematics and Applied Statistics, University of Wollongong\\
Wollongong, Australia
}

\advance\topmargin by -\headheight
\advance\topmargin by -\headsep
\evensidemargin=\oddsidemargin
\makeatletter
\def\section{\@startsection{section}{1}{\z@}{-3.5ex plus -1ex minus
  -.2ex}{2.3ex plus .2ex}{\large\bf}}
\def\subsection{\@startsection{subsection}{2}{\z@}{-3.25ex plus -1ex
  minus -.2ex}{1.5ex plus .2ex}{\normalsize\bf}}
\makeatother

\numberwithin{equation}{section} 

\theoremstyle{plain} 

\theoremstyle{definition} 

\newtheorem*{ass*}{Standing assumption}
\newtheorem*{example*}{Example} 
\theoremstyle{remark} 

\newcommand{\C}{\mathbb{C}}   
\newcommand{\stroke}{\mathbin|}     

\def\pairL_#1(#2|#3){{}_{#1}(#2\stroke#3)} 
\def\pairR(#1|#2)_#3{(#1\stroke#2)_{#3}} 
\def\scal<#1|#2>{\langle#1\stroke#2\rangle} 

\newbox\ncintdbox \newbox\ncinttbox 
	\setbox0=\hbox{$-$}
	\setbox2=\hbox{$\displaystyle\int$}
	\setbox\ncintdbox=\hbox{\rlap{\hbox
		to \wd2{\hskip-.125em \box2\relax\hfil}}\box0\kern.1em}
	\setbox0=\hbox{$\vcenter{\hrule width 4pt}$}
	\setbox2=\hbox{$\textstyle\int$}
	\setbox\ncinttbox=\hbox{\rlap{\hbox
		to \wd2{\hskip-.175em \box2\relax\hfil}}\box0\kern.1em}

\begin{document}

\maketitle

\vspace{-2pc}

\begin{abstract}
We show by (counter)example that the intersection of complemented submodules in a Hilbert $C^*$-module is not necessarily complemented, 
answering an open question from \cite{MR}.
\end{abstract}



\parskip=6pt
\parindent=0pt
\allowdisplaybreaks
\section*{The intersection of complemented submodules}
Let $X$ be a right Hilbert $C^{*}$-module over the $C^{*}$-algebra $B$ and $M\subset X$ a closed $B$-submodule. The \emph{orthogonal complement} of $M$ is the submodule
\[M^{\perp}:=\left\{x\in X : \forall m\in M\,\,\langle x,m\rangle =0\right\}.\] Recall that $M$ is \emph{complemented} if the map 
\[M\oplus M^{\perp}\to X,\quad (m,x)\mapsto m+x,\]
is an isomorphism of Hilbert $C^{*}$-modules. Given two complemented submodules $M,N\subset X$, their intersection $M\cap N\subset X$ is a closed $B$-submodule.
%
In \cite{MR} the second and third named authors posed: 

``To our knowledge, it is an open question whether the intersection of complemented submodules is again complemented.'' 

In the present note we settle this question in the negative. The counterexample below was provided by the first named author. 

\begin{example*}
Consider the Hilbert $C^*$-module $X := C([-\pi/2,\pi/2],\C^2)$ over the $C^*$-algebra $B := C([-\pi/2,\pi/2])$. 
For $t\in [-\pi/2,\pi/2]$ define projection-valued functions
\[
P(t)=\begin{pmatrix}1&0\\0&0\end{pmatrix}
\quad\mbox{and}\quad
Q(t)=\left\{\begin{array}{ll}\begin{pmatrix} \cos^{2}t & \sin t \cos t \\ \sin t \cos t &\sin^{2}t\end{pmatrix}&t\geq0\\\begin{pmatrix}1&0\\0&0\end{pmatrix}&t<0\end{array}\right..
\]
For the complemented submodules $M:={\rm Im}(P)$ and $N:={\rm Im}(Q)$ in $X$ we have
\[
M\cap N=\left\{\begin{pmatrix} x \\ 0\end{pmatrix}:
x \in C([-\pi/2,\pi/2]),\ 
x(t)=0\ \mbox{for }t\geq 0\right\} ,
\]
where the behaviour at $t=0$ comes from continuity of any $x$ in the images of $P,Q$. So then continuity again says that the complement is
\[
(M\cap N)^\perp=\left\{\begin{pmatrix} y \\ z\end{pmatrix}:
y,z \in C([-\pi/2,\pi/2]),\ 
y(t)=0\ \mbox{for }t\leq 0\right\}.
\]
Hence any vector $\begin{pmatrix} y \\z \end{pmatrix}\in C([-\pi/2,\pi/2],\C^2)$ with $y(0)\neq0$
is not in $(M\cap N) + (M\cap N)^\perp$.
Thus $M\cap N$ is a closed submodule which is not complemented. 
\end{example*}


\begin{thebibliography}{9999}
%
%
%
%
%
%
%
%
%
%
%
%
%
\bibitem[MR]{MR} B. Mesland, A. Rennie, {\em The Friedrichs angle and alternating projections in Hilbert $C^{*}$-modules}, J. Math. Anal. Appl., {\bf 516} (2022), 126474.

%
%
%
%
%

\end{thebibliography}
\end{document}